\documentstyle{amsart}
\input{bull-art}
\newtheorem{thm}{Main Theorem} 
\newtheorem{theorem}[subsection]{Theorem}
\newtheorem{conjecture}[subsection]{Conjecture}
\newtheorem{hypo}[subsection]{Hypothesis}

\theoremstyle{definition}
\newtheorem{maximal}{Maximal parabolic geometries} 

\newtheorem{graphs}{Graphs} 
\newtheorem{unique}{Uniqueness of certain sporadics} 

\newtheorem{generators}{Generators and relations}

\newcommand{\CHAR}{\operatorname{char}}
\newcommand{\SP}{\operatorname{Sp}}
\newcommand{\G}{{\cal G}}
\newcommand{\A}{{\cal A}}
\newcommand{\Aut}{\operatorname{Aut\,}}
\newcommand{\res}{\operatorname{res}}
\newcommand{\wreath}{\operatorname{wr}}

\mathchardef\tnode="020E 

\def\arc{
  \hbox{\kern -0.15em
  \vbox{\hrule width 2.5em height 0.6ex depth -0.5 ex}
  \kern -0.33em}}

\def\darc{
  \rlap{\lower0.2ex\arc}{\raise0.2ex\arc}}

\def\stroke#1{
  \kern 0.05em
  \rlap\arc{{\textstyle{#1}}\atop\phantom\arc}
  \kern -0.22em}

\def\dstroke#1{
  \kern 0.05em
  \rlap\darc{{\textstyle{#1}}\atop\phantom\darc}
  \kern -0.22em}

\def\centerscript#1{
  \setbox0=\hbox{$\tnode$}
  \hbox to \wd0{\hss$\scriptstyle{#1}$\hss}}

\def\node{
  \def\super{}
  \def\sub{}
  \futurelet\next\dolabellednode}

  \let\sp=^
  \let\sb=_

  \def\dolabellednode{%
    \ifx\next\sb\let\next\getsub
    \else
      \ifx\next\sp\let\next\getsuper
      \else\let\next\donode
      \fi
    \fi
    \next}

  \def\getsub_#1{\def\sub{#1}\futurelet\next\dolabellednode}
  \def\getsuper^#1{\def\super{#1}\futurelet\next%
\dolabellednode}

  \def\donode{%
    \rlap{$\mathop{\phantom\tnode}\limits_{\centerscript{%
\sub}}
    ^{\centerscript{\super}}$}\tnode}

\def\varcdn{
  \kern -0.03em\vbox{\kern -0.5ex
  \hbox to \wd0{\hss\vrule width 0.04em depth 5.8ex\hss}
  \kern -0.3ex
  \hbox{$\tnode$}}}

\begin{document}
\def\currentvolume{31}
\def\currentissue{2}
\def\currentyear{1994}
\def\currentmonth{October}
\def\copyrightyear{1994}
\def\currentpages{173-184}
\ratitle
\title[The flag-transitive tilde]{The flag-transitive tilde and Petersen-type 
geometries are all known}
\author{A. A. Ivanov}
\address{Institute for Systems Analysis,
Russian Academy of Sciences, 9,
Prospect 60 Let Oktyabrya, 117312, Moscow,
Russia}
\email[A. Ivanov]{ivanov@@cs.vniisi.msk.su\newline\indent
    {\defaultfont{\it E-mail address}, S. Shpectorov: 
}ssh@@cs.vniisi.msk.su}
\author{S. V. Shpectorov}
\date{March 11, 1993}

\maketitle
\begin{abstract}
We announce the classification of two related classes of 
flag-transitive
geometries. There is an infinite family of such 
geometries, related to the
nonsplit extensions $3^{[{n\atop 2}]_{_2}}\cdot 
\SP_{2n}(2)$, and twelve sporadic
examples coming from the simple groups $M_{22}$, $M_{23}$, 
$M_{24}$, $He$,
$Co_1$, $Co_2$, $J_4$, $BM$, $M$ and the nonsplit extensions
$3\cdot M_{22}$, $3^{23}\cdot Co_2$, and $3^{4371}\cdot BM$.
\end{abstract}

\section{Introduction}
\par
An {\em incidence system} $\G$ is a set
of elements, each having a particular type (or coloring), 
on which a
reflexive, symmetric, binary incidence relation is 
defined. Two different 
elements of the same type can never be incident. A {\em 
flag}
is a set of pairwise incident elements. The {\em residue} 
of a flag $F$, 
denoted $\res(F)$, is the incidence system naturally 
formed by all the 
elements not in $F$, which are incident to all the 
elements from $F$. 
The rank of $\G$ is the number of different types 
presented in $\G$.
An incidence system $\G$ is {\em connected} if the graph 
on $\G$ defined by 
the incidence is connected. 

A {\em geometry} is an incidence system for which two 
further properties hold:
(1) every maximal flag contains elements of all types; (2) 
every residue of 
rank at least 2 (including the geometry itself) is 
connected. The latter 
condition is known as {\em strong} or {\em residual 
connectedness}. Both 
these properties are inductive, i.e., the residues of 
geometries are 
themselves geometries. This sort of inductiveness was used 
in the definition 
of the diagram of a geometry. An example of
such a diagram can be found below. 
The vertices of a diagram are all the types of the 
corresponding geometry 
(geometries). If you remove from the diagram all the types 
presented 
in a flag $F$ and all the edges on these types, then the 
remaining 
diagram will describe $\res(F)$\<; in particular, the 
edges of the diagram 
describe the rank 2 residues of the geometry. As the type 
set we usually 
take $\{1,2,\ldots,n\}$, 
where $\bold n$ is the rank and
the types in the diagram increase 
from left to right. For more 
information about the axioms and the properties of the 
geometries and diagrams see \cite{bue1,tit2}.

A {\em flag-transitive} automorphism group $G$ is a group 
of permutations on 
a geometry $\G$, which preserves types and the incidence, 
and, moreover, acts 
transitively on the set of maximal flags of $\G$. In this 
situation we also
call $\G$ flag-transitive. We often say that $\G$ is the 
geometry of the 
group $G$ and denote it $\G=\G(G)$. 

The systematic study of the geometries of sporadic simple 
groups
was initiated in \cite{bue1}. The idea was to develop a 
geometric theory of
sporadic groups close to Tits's theory of buildings for 
the Chevalley groups
\cite{tit1,tit2}. In the latter case the diagrams of the 
buildings 
naturally correspond to the Dynkin diagrams of the 
corresponding 
Chevalley groups, and just like the Dynkin diagrams, the 
geometric 
diagrams of the buildings contain enough information to 
reconstruct the 
buildings and classify them. 

By now, examples of flag-transitive geometries are known 
for all twenty-six sporadic 
groups. The diagrams of these geometries contain many new 
edges. Some (but not many) 
of the sporadic geometries were characterized by their 
diagrams, usually under the 
flag-transitivity assumption (cf. \cite{pas2,bue2}).

In this paper we announce the classification of the 
flag-transitive
geometries belonging to the following two types of diagrams:
\medskip
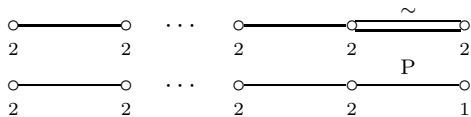
\setlength{\unitlength}{1mm}
$$
\begin{picture}(60,5)(0,-2.5)
%
%
\put(0,0){\makebox(0,0){$\circ$}}
\put(0,-3){\makebox(0,0){$\scriptstyle{2}$}}
\put(0.7,0.1){\line(1,0){13.5}}
\put(15,0){\makebox(0,0){$\circ$}}
\put(15,-3){\makebox(0,0){$\scriptstyle{2}$}}
\put(22.5,0){\makebox(0,0){$\cdots$}}
\put(30,0){\makebox(0,0){$\circ$}}%
\put(30,-3){\makebox(0,0){$\scriptstyle{2}$}}
\put(30.7,0.1){\line(1,0){13.5}}
\put(45,0){\makebox(0,0){$\circ$}}
\put(45,-3){\makebox(0,0){$\scriptstyle{2}$}}
\put(45.5,-0.6){\line(1,0){14}}
\put(45.5,0.7){\line(1,0){14}}
\put(52.5,2){\makebox(0,0){$\scriptstyle{\sim}$}}
\put(60,0){\makebox(0,0){$\circ$}}
\put(60,-3){\makebox(0,0){$\scriptstyle{2}$}}
\end{picture}
$$
$$
\begin{picture}(60,5)(0,-2.5)
%
%
\put(0,0){\makebox(0,0){$\circ$}}
\put(0,-3){\makebox(0,0){$\scriptstyle{2}$}}
\put(0.7,0.1){\line(1,0){13.5}}
\put(15,0){\makebox(0,0){$\circ$}}
\put(15,-3){\makebox(0,0){$\scriptstyle{2}$}}
\put(22.5,0){\makebox(0,0){$\cdots$}}
\put(30,0){\makebox(0,0){$\circ$}}
\put(30,-3){\makebox(0,0){$\scriptstyle{2}$}}
\put(30.7,0.1){\line(1,0){13.5}}
\put(45,0){\makebox(0,0){$\circ$}}
\put(45,-3){\makebox(0,0){$\scriptstyle{2}$}}
\put(45.7,0.1){\line(1,0){13.5}}
\put(52.5,2.5){\makebox(0,0){{\small\rm P}}}
\put(60,0){\makebox(0,0){$\circ$}}
\put(60,-3){\makebox(0,0){$\scriptstyle{1}$}}
\end{picture}
$$

\bigskip
\noindent
where
{$
\begin{picture}(16,3)(0,-0.8)
\put(0,0){\makebox(0,0){$\circ$}}
\put(15,0){\makebox(0,0){$\circ$}}
\put(0.4,-0.6){\line(1,0){14.1}}
\put(0.4,0.7){\line(1,0){14.1}}
\put(7.5,2){\makebox(0,0){$\scriptstyle{\sim}$}}
\end{picture}
$}
(resp.
{$
\begin{picture}(16,3)(0,-0.8)
\put(0,0){\makebox(0,0){$\circ$}}
\put(15,0){\makebox(0,0){$\circ$}}
\put(0.7,0.1){\line(1,0){13.5}}
\put(7.5,2.5){\makebox(0,0){{\small\rm P}}}
\end{picture}
$}
)
denotes the triple cover of the generalized quadrangle for 
$\SP_4(2)$
(resp. the geometry of edges and vertices of the Petersen 
graph).
We  call these geometries, respectively, {\em tilde} (for 
short, $T$-) and 
{\em Petersen-type} (for short, $P$-) geometries.

We present the list of $T$- and $P$-geometries and some of 
their properties
in Table 1. In the fourth column we indicate the residue 
of an element of type 1; 
this defines a tree structure on the set of $T$- and 
$P$-geometries. 
Other entries will be explained later.

{\fontsize{9}{11pt}\selectfont 
\begin{table}[t]\caption{$T$- and $P$-geometries}
\medskip
\centering
\begin{tabular}{|c|c|c|c|c|c|c|}
\hline
\# & rank & $\Aut{\cal G}$                         & Res. 
& 2-cov. & Subgeom.     
&Nat. Reps.         \\
\hline
\hline
$T_0$ & 2    & $3\cdot \text{Sp}_4(2)$                     
   & $-$  & $-$    & $P_0$            & 6+5                
 \\
$T_0(n)$ & $n>2$  & $3^{[{n\atop 2}]_{_2}}\cdot 
\text{Sp}_{2n}(2)$ & $T_0(n-1)$   & $-$    & $-$           
& $2^n(2^n-1)+$        \\
&&&&&& $(2n+1)$        \\
$T_1$ & 3    & $M_{24}$                                & 
$T_0$   & $-$    & $P_1$            & 11                   
\\
$T_2$ & 3    & $He$                                    & 
$T_0$   & $-$    & $-$           & 52                    
\\
$T_3$ & 4    & $Co_1$                                  & 
$T_1$   & $-$    & $P_4$            & 24                   
 \\
$T_4$ & 5    & $M$                                     & 
$T_3$   & $-$    & $P_7$            & $-$               \\
\hline
\hline
$P_0$ & 2    & $S_5$                                   & 
$-$  & $-$    &               & 6            \\
$P_1$ & 3    & $\Aut M_{22}$                           & 
$P_0$   & $P_2$     &
$\text{Sp}_4(2)$     & 11            \\
$P_2$ & 3    & $3\cdot \Aut M_{22}$                    & 
$P_0$   & $-$    & $T_0$            & 12+11             \\
$P_3$ & 4    & $M_{23}$                                & 
$P_1$   & $-$    & $A_7$         & $-$            \\
$P_4$ & 4    & $Co_2$                                  & 
$P_1$   & $P_5$     &
$\text{Sp}_6(2)$     & 23            \\
$P_5$ & 4    & $3^{23}\cdot Co_2$                      & 
$P_2$   & $-$    & $T_0(3)$         & 0+23   \\
$P_6$ & 4    & $J_4$                                   & 
$P_2$   & $-$    & $T_1$            & $-$           \\
$P_7$ & 5    & $BM$                                    & 
$P_4$   & $P_8$     &
$\text{Sp}_8(2)$     & $-$           \\
$P_8$ & 5    & $3^{4371}\cdot BM$                      & 
$P_5$   & $-$    & $T_0(4)$         & $-$           \\ 
\hline
\end{tabular}
\end{table}
}%

\begin{thm}
Every flag-transitive T- or P-geometry is isomorphic
to one of the finite geometries listed in Table \rom{1}.
\end{thm}

For a flag-transitive group, by analogy with the case of 
buildings, the 
stabilizers of nonempty flags are called the {\em 
parabolic} subgroups. 
The rank of a parabolic is the rank of the residue of the 
corresponding flag.
Let $F$ be a maximal flag, and let $M_1,M_2,\ldots,M_n$ be 
the stabilizers 
of the elements from $F$ ($M_i$ is the stabilizer of the 
element of type $i$). 
These stabilizers form the {\em amalgam of maximal 
parabolic subgroups}. 
By {\em amalgam} we understand just a collection of groups 
with common identity
element and with the group operations consistent on 
intersections.
The amalgam of maximal parabolics defines all the proper 
residues of 
the geometry, and, in particular, it defines the diagram. 
The flag-transitive 
geometries giving rise to the same amalgam are connected 
with each other 
by the operations of taking covers/quotients. A {\em 
covering} is a surjective 
morphism which is an isomorphism on every proper residue. 
Similarly, 
an {\em $s$-covering} is a surjective morphism which is an 
isomorphism on 
every residue up to rank $s$. There is a number of 
nontrivial 2-coverings 
between our geometries; they are shown in the fifth column 
of Table 1.

Within our approach the proof of the Main Theorem can be 
formally divided into 
two steps.  First we classify possible amalgams of maximal 
parabolic subgroups (local characterization). The 
principal result is that 
every such amalgam comes from a known example.  Second we 
determine the geometries/groups corresponding to the known 
amalgams. 
It amounts to proving that the final list of geometries is 
closed 
with respect to
 taking flag-transitive covers/quotients. It turned out 
that the original 
list of geometries was incomplete, and so we also faced a 
problem of constructing 
new $T$- and $P$-geometries related to strange nonsplit 
extensions of 
some sporadic and serial simple groups.

\section{The geometries}
\par
The $T$-geometries $\G(M_{24})$, $\G(He)$, $\G(Co_1)$, and 
$\G(M)$ 
were constructed in \cite{rst} as minimal parabolic 
geometries for the corresponding groups
using their maximal parabolic geometries from \cite{rsm1}. 
The $P$-geometries
$\G(M_{22})$, $\G(3 \cdot M_{22})$, $\G(M_{23})$, 
$\G(Co_2)$, and $\G(BM)$ were
constructed in \cite{ivn1} in terms of graphs. The 
truncations of these geometries 
by the elements of maximal  type
coincide with the minimal parabolic geometries from 
\cite{rst}.

Suppose $G<\Aut H$, where $H$ is a group, and suppose $H$ 
contains a subgroup $E\cong 2^n$, such 
that $N_G(E)$ induces on $E$ the whole group $L_n(2)$. 
Define an incidence system $\G=\G(G,H,E)$ 
whose elements are all the subgroups of $H$ conjugate 
under $G$ to nonidentity subgroups of $E$;
the type is equal to the 2-rank, and the incidence is 
defined by inclusion. Then $\G$ fulfills the 
axioms of a geometry (except for the strong connectedness, 
which must be checked separately), 
$G$ acts flag-transitively on $\G$, and $\G$ belongs to a 
string diagram (i.e., without 
branches or loops) in which the residue of 
an element of type $n$ is isomorphic to the projective 
geometry of the proper subgroups of $E$.

The  following configurations give rise to $T$- and 
$P$-geometries; 
cf. \cite{is2}. In the first 
class of examples $H$ is a suitable $GF(2)G$-module. For 
$G\cong 3\cdot \SP_4(2)$ or $S_5$ it is 
the natural module of $\Gamma L_3(4)>G$ (the hexacode); 
for $G\cong M_{24}$ or $\Aut(M_{22})$ it is 
the Golay cocode; for $G\cong 3\cdot\Aut(M_{22})$ it is 
the natural module of $\Gamma U_6(2)>G$; 
for $G\cong Co_1$ or $Co_2$ it is $\Lambda/2\Lambda$, 
where $\Lambda$ is the Leech lattice; 
finally, for $G \cong He$ it is a 51-dimensional rational 
module, taken modulo 2 \cite{ms}.
Notice that the conjugates of $E$ need not span the whole 
of $H$, say; in the case of $G\cong Co_2$ 
the conjugates of $E$ span a submodule of codimension 1 in 
$\Lambda/2\Lambda$.
In its turn, this 
submodule has a 22-dimensional quotient, which can also be 
taken as $H$.

For $G\cong J_4$, $BM$, or $M$ the role of $H$ is played 
by $G$ itself with the natural action by 
conjugation. In each of these cases $G$ contains an 
involution $\langle\tau\rangle$, such that 
$C=C_G(\tau)$ has the form $C=Q.A$ where $Q$ is an 
extraspecial group of order $2^{1+m}$ 
for $m=12$, 22, or 24, respectively. Let $\bar 
Q=Q/\langle\tau\rangle$. Then the action of $A$ 
(isomorphic, respectively, to $3\cdot\Aut(M_{22})$, 
$Co_2$, or $Co_1$) on $\bar Q$ corresponds to 
a certain configuration from the previous paragraph. In 
particular, there is a subgroup $\bar E$
in $\bar Q$, which gives rise to the geometry $\G(A)$. We 
take $E$ to be the full preimage of 
$\bar E$ in $Q$. In each case it is easy to see that 
$N_G(E)$ induces on $E$ the full linear 
group, and it is almost by the definition that the residue 
in $\G(G,G,E)$ of the element 
$\langle\tau\rangle$ coincides with $\G(A,\bar Q,\bar E)$.

An additional class of examples can be constructed as 
$i$-covers. Let $G\cong
\SP_{2n}(2)$, 
$\Aut(M_{22})$, $Co_2$, or $BM$ and $\G=\G(G)$ be the 
corresponding geometry.
By $\G(\SP_{2n}(2))$ 
we mean the classical $C_n(2)$-geometry. Consider in $G$ a 
subgroup $L$ isomorphic to $\Omega^-_{2n}(2).2$, 
$L_3(4).2$, $U_6(2).2$, or $2\cdot^2\!E_6(2).2$, 
respectively. 
Then $L/L' \cong Z_2$
(here $L'$ is the derived group of $L$), and there is a 
unique nontrivial
1-dimensional $GF(3)L$-module $X$ whose kernel is $L'$. 
Let $Y$ be the
$GF(3)G$-module induced from $X$ and $\hat G \cong Y.G$
be a certain extension of $G$ by $Y$. We consider
the split extension in the first case and
nonsplit extensions in the other cases. For $L' \cong 
L_3(4)$, $U_6(2)$, or
$2 \cdot ^2E_6(2)$ the 3-part of the Schur multiplier of 
$L'$ is of order 3.
Moreover, every element from $L - L'$ when extended to an 
automorphism of
the nonsplit triple cover of $L'$ inverts the center. Now by
the Eckmann-Shapiro lemma
(cf. Shapiro's lemma in \cite{bro}) in each of these three 
cases there
exists a unique nonsplit extension $\hat G$.
It can be shown that the amalgam of the rank $i$ parabolic 
subgroups from $G$ ($i=1$ for $G\cong \SP_{2n}(2)$, and 
$i=2$ otherwise) is embedded in $\hat G$. 
In the first case there is a trivial embedding into a 
complement to $Y$, and we consider another 
one. The embedded amalgam generates in $\hat G$ a nonsplit 
extension 
$3^{[{n\atop 2}]_{_2}}\cdot \SP_{2n}(2)$, 
$3\cdot\Aut(M_{22})$, $3^{23}\cdot
Co_2$,
 or $3^{4371}\cdot BM$, 
respectively (where $[{n\atop 2}]_{_2}=(2^n-1)(2^n-2)/6$). 
This extension corresponds to a 1-covering 
of $\G(\SP_{2n}(2))$ by a flag-transitive $T$-geometry in 
the first case and to a 2-covering of $\G$ by a 
$P$-geometry in the other three cases. These constructions 
were accomplished in \cite{is6,shp3,is7}. 
The symplectic series of $T$-geometries was also 
independently constructed by U. Meierfrankenfeld 
\cite{mei}. Notice that the exceptional $C_3(2)$-geometry 
for the group $A_7$ does not 
have a 1-cover which is a $T$-geometry, as was checked in 
\cite{is3}.

As indicated in column six of Table 1 on page 175
many of the $T$- and $P$-geometries
contain subgeometries of $T$-, $P$-, or $C_n(2)$-type. These
subgeometries play a crucial inductive role in the 
classification.

\section{Simple connectedness}
\par
Within the approach we choose, the
following result constitutes a very important step of the 
classification.
\begin{theorem}
The set of geometries in Table $1$ with rank at least
$3$ is closed with respect to taking flag-transitive 
covers and quotients.
\end{theorem}

Of course, the hard part of this statement is that every 
geometry in the table with 
rank at least 3, except $\G(M_{22})$, is simply connected 
(i.e., has no nontrivial 
coverings). This is equivalent to the fact that the 
automorphism group of the
geometry coincides with the universal completion $U(\A)$ 
of the
corresponding amalgam $\A$ of maximal parabolic subgroups.

Since $U(\A)$ can be defined in terms of generators and 
relations, 
for its identification one can use the coset enumeration 
algorithm.
This was implemented in \cite{hei} to check the simple 
connectedness of $\G(M_{24})$, 
$\G(He)$, $\G(3^7 \cdot \SP_6(2))$. Later an independent 
computer-free proof for the case 
of $\G(M_{24})$ was
found by the first author.

Another strategy goes back to \cite{ron1} and relies on 
analysis of cycles
in the {\em collinearity} graph $\Gamma$ of $\G$. The 
collinearity graph has $\G^1$ as 
the set of vertices with two vertices adjacent if they are 
incident to a common 
element from $\G^2$. In many important cases it can be 
shown that a covering of 
$\G$ induces a covering of its collinearity graph  and 
that with respect 
to the induced covering all triangles are null-homotopic. 
In this case to prove the
simple connectedness of $\G$, it is sufficient to show 
that $\Gamma$ is {\em triangulable}, 
which means by definition that every cycle in $\Gamma$ can 
be decomposed into a product of 
triangles. 

Proceeding by induction on the
rank and the number of elements,
 we can assume that all $T$- and $P$-subgeometries
(as well as $C_n(2)$-subgeometries, if any) in $\G$ are 
simply connected.
Then each cycle of $\Gamma$  which completely lies in a 
subgeometry is 
null-homotopic,
 and the analysis can be simplified considerably. This 
scheme
was realized in \cite{is3} for $\G(M_{23})$ and in 
\cite{shp3} for $\G(Co_2)$.

For larger geometries it turned out to be more convenient 
to consider a different 
graph $\Sigma=\Sigma(\G)$ (the intersection graph of 
subgeometries).
The vertices of $\Sigma$ are subgeometries, and two 
subgeometries are adjacent if 
they have ``large'' intersection.
It was shown that every covering of $\G$ induces a 
covering of $\Sigma$ and that 
all triangles are null-homotopic with respect to the 
induced covering. After that 
it was proved that $\Sigma$ is
triangulable and hence $\G$ is simply connected. 
This approach was realized for $\G(Co_1)$, $\G(J_4)$,
$\G(BM)$, and $\G(M)$ in \cite{ivn9,ivn7,ivn8,ivn6}, 
respectively.

For the remaining geometries the simple connectedness 
proof involves both
combinatorial and group-theoretical arguments. In these 
cases
$O_3(G) \ne 1$, and the homomorphism $\phi: G \to \bar 
G=G/O_3(G)$ induces a morphism
of $\G(G)$ onto a geometry $\G(\bar G)$, and the latter 
one is already known
to be simply connected. We prove that for the amalgam $\A$ 
of maximal
parabolics the inequality $|U(\A)| \le |G|$ holds 
\cite{is6,is7,shp3}.
\section{Natural representations}
\par
The geometries $\G(G,H,E)$ from Section 2, where  
$H$ is a $GF(2)G$-module, are constructed in their 
``natural representations''.
By a {\em natural representation} of a $T$-, $P$-, or 
$C_n(2)$-geometry $\G$
we understand its morphism $\phi$ into the projective 
geometry of the nontrivial subspaces of
a $GF(2)$-space, such that for $x\in\G^i$ we have 
$\dim(\phi(x))=i$ and for
$j\le i$ the restriction of $\phi$ to the set of $j$-type 
elements from $\res_\G(x)$
is a bijection onto the set of $j$-dimensional subspaces 
of $\phi(x)$. 

As all our geometries  are of $GF(2)$-type (three points 
on a line), the natural 
representations of a fixed geometry $\G$ (if they exist at 
all) are 
controlled by a particular one known as the {\em universal 
natural representation}. 
The corresponding $GF(2)$-vector space can be defined as 
the largest space 
spanned by vectors $v_p$, $p\in\G^1$, such that $v_a+v_b+
v_c=0$ whenever 
$\{a,b,c\}=\res_\G(\ell)^1$ for an element $\ell\in\G^2$. 
For $\G=\G(G)$ this largest vector space $UM(\G)$ can be 
considered as a $G$-module and 
 is called the {\em universal natural module}.

The above notions proved to be useful in the local 
characterization of the amalgams of maximal
parabolics, where we need information on certain natural 
modules. 
In principle, the determination of the dimension of 
$UM(\G)$ is a problem of linear algebra 
over $GF(2)$. In particular, for the groups $G\cong 3\cdot 
S_6$ and $S_5$ we can easily determine 
$\dim(UM(\G(G)))$.  Another and much more impressive 
computation is due to B.~McKay \cite{mck}, who
established $\dim(UM(\G(He)))=52$ by considering on 
computer a system of 437,325 equations 
over 29,155 variables.

In all the other cases the proofs are computer-free 
\cite{is2,rsm3,smi,is4,is5,shp4}. 
Summarizing these results, we obtain the following.
\begin{theorem}
The dimensions of the universal natural representations
of $T$- and $P$-geometries are as given in column seven of 
Table $1$.
\end{theorem}

If, for the automorphism group $G$, we have $O_3(G) \ne 
1$, then in Table 1 the dimension of 
the universal natural module is given in two summands: the 
commutant and the  centralizer 
of $O_3(G)$, respectively.

Of particular interest are the arguments for the groups 
$G\cong J_4$, $BM$,
$3^{4371}\cdot BM$, and $M$. If in the above definition of 
$UM(\G)$ we do not assume 
that the elements $v_p$ commute, we obtain the definition 
of the {\em universal natural group} 
$UG(\G)$. If $G$ is one of the above groups and 
$\G=\G(G)$, then $G$ is a quotient of 
$UG(\G)$ (cf. Section 2). Let $\xi:UG(\G)\to G$ be the 
corresponding homomorphism. 
It was shown that the restriction of $\xi$ to a certain 
subgroup $Q\le UG(\G)$ 
is an isomorphism and $Q'$ contains $v_p$ for some 
$p\in\G^1$. This readily implies 
that $UG(\G)$ has no abelian quotients, i.e., $UM(\G)$ is 
trivial. 
The nonexistence of the natural representations of $\G(M)$ 
answers a question posed in \cite{str}.

\begin{conjecture}
Let $G\cong J_4$, $BM$, or $M$. Then $UG(\G(G))$ is 
isomorphic to 
$J_4$, $2\cdot BM$, or $M$.
\end{conjecture}

At the moment the conjecture is 
proved for $G\cong BM$ and $M$ \cite{ips}.

\section{Local characterization}
\par
On the stage of local characterization we prove that the 
amalgam 
$\A=\{M_1,...,M_n\}$ of maximal parabolics corresponding 
to a flag-transitive action on a $T$- or 
$P$-geometry is isomorphic to that from a relevant known 
example.

Partial results on local characterization of 
$T$-geometries were proved in
\cite{hei,row1,row2,row3,tim}. Compared with these papers, 
we use a somewhat
different approach; 
in particular, we do not assume finiteness of parabolics. 
The approach was developed in 
\cite{shp1,shp2} for $P$-geometries, and it works for 
$T$-geometries as well.
To simplify the notation, 
let us consider only one of the two types of geometries, 
say, $T$-geometries.

Let $\Delta$ be the graph on $\G^n$ ($n$ is the rank of 
$\G$), in which two elements are 
adjacent exactly when they are incident to a common 
element of type $n-1$ ($\Delta$ is called the 
{\em derived
graph} of $\G$). For $y\in\G^i$, $i<n$, let $\Sigma(y)$ be 
the subgraph of $\Delta$ induced by the 
elements (vertices) incident to $y$. Let ${\cal O}$ be the 
set of the subgraphs $\Sigma(y)$. 
For $y$ of type $i<n-1$ the subgraph $\Sigma(y)$ 
is naturally isomorphic to the derived graph of the 
residual $T$-geometry of rank $n-i$, related to 
$y$. If $y\in\G^{n-1}$, then $\Sigma(y)$ is a 3-clique. It 
follows from the diagram  of $\G$ that the 
subgraphs $\Sigma\in{\cal O}$ containing a particular 
vertex $a$ (we denote this set of subgraphs 
by ${\cal O}_a$) naturally form an $(n-1)$-dimensional 
projective space over $GF(2)$. The graph 
$\Delta$ together with the set ${\cal O}$ gives another 
realization of $\G$.

The group $G$ acts naturally on $\Delta$ preserving ${\cal 
O}$ and having $M=M_n$ as the stabilizer 
of a vertex $a$. Let $K_s$ be the subgroup of $M_n$ fixing 
all the vertices at
a distance at most $s$ 
from $a$. Consider the following condition.

\vspace{0.2cm}

\noindent
$(*)$ $K_{n-1}$ (where, as above, $n$ is the rank of $\G$) 
is a group of order at most 2.

\vspace{0.2cm}

When $(*)$ holds for all the residual $T$-geometries up to 
rank $i<n$, we can 
prove that $K_i/K_{i+1}$ is an irreducible $GF(2)$-module 
for the quotient $L_n(2)$ of 
$M$ (this quotient is the action of $M$ on ${\cal O}_a$), 
of dimension 0, 1, or ${n\choose i}$. 
The proof goes as follows. Every vertex at distance $i$ 
from $a$ is covered by a subgraph 
$\Sigma(y)\in{\cal O}_a$ for a $y$ of type $n-i$. By 
$(*)$, $K_i$ induces on $\Sigma(y)$ 
the action of order at most 2. If the action is trivial, 
$K_i/K_{i+1}$ is trivial as well. 
Otherwise, in the dual of $K_i/K_{i+1}$ there is an orbit 
indexed by the $y$\<'s, which are simply 
all the elements of type $i$ in the projective space 
$\res_\G(a)\cong{\cal O}_a$. We can check 
that the vectors from this orbit (if not all equal) 
possess some 
natural 3-term linear relations, which immediately leads 
to the identification of $K_i/K_{i+1}$.

If all the proper residual $T$-geometries have the 
property $(*)$, then some further simple arguments 
show $K_{n+1}=1$ and enable us to determine the structure 
of $K_n$. As soon as all the chief factors of 
$M$ are known, reconstruction of the possible amalgams 
$\A$ is only a matter of technique.

Since the condition $(*)$ holds for all the known 
$T$-geometries, except for the 
terminal (in the tree of $T$-geometries) geometry $\G(M)$, 
we can inductively 
reconstruct all the amalgams corresponding to the existing 
$T$-geometries. 
It remains to prove that the $T$-geometry $\G(M)$ has no 
further extension. 
For such an extension we can determine, as above, 
four first sections
$K_i/K_{i+1}$. However, 
we cannot immediately conclude that $K_{n+1}=1$ (here 
$n=6$). Let us consider the parabolic $M_1$ which is 
the stabilizer in $G$ of an element $u\in\G^1$. Let $N$ be 
the kernel of $M_1$ acting on $\res_{\G}(u)$ 
(equivalently, on $\Sigma(u)$) and $N_1$ be the kernel of 
$N$ acting on the set of elements from 
$\G^1$ which are incident with $u$ to a common element 
from $\G^2$. Then $M_1/N$ is a flag-transitive 
action on $\res_{\G}(u)$ and is known to be isomorphic to 
the Monster group $M$ by the
inductive hypothesis. We prove that $N/N_1$ is a 
$GF(2)$-space and its dual realizes a 
natural representation of $\res(u)$. By Theorem 4.1 this 
module must be trivial, i.e., $N=N_1$. 
Then we prove, using $(*)$ for the residues up to rank 4, 
that $N$ stabilizes all the vertices 
at distance $\le 4$ from $\Sigma(u)$, which means that 
$N\le K_4$. However, the 
2-part of the order of $M_n/K_4$ is greater than that of 
$M_1/N\cong M$, unless the section 
$K_3/K_4$ is already trivial or 1-dimensional. In the 
latter case we obtain that $K_4=1$ and 
eventually establish a contradiction.

\section{Applications}
Let us discuss some consequences of the Main Theorem and 
its proof.

\begin{maximal}
The simple connectedness of the minimal para\-bolic
geometries of the
groups $M_{24}$, $He$, $Co_1$, $Co_2$, $J_4$, $BM$, and 
$M$ implies the
simple connectedness of the corresponding maximal 
parabolic geometries
(proved in \cite{ron2} and \cite{seg1} for the cases 
$M_{24}$ and $Co_1$,
respectively). The fact that the universal covering of 
$\G(M_{22})$ is
isomorphic to $\G(3 \cdot M_{22})$ implies the simple 
connectedness of
the maximal parabolic geometry for $M_{22}$.
\end{maximal}

\begin{graphs}
First examples of $P$-geometries were constructed in terms 
of
the corresponding derived graphs \cite{ivn1}. These graphs 
satisfy the
following.
\end{graphs}

\begin{hypo}
$\Gamma$ is a graph whose girth  \rom{(}the length of
the shortest cycle\rom{)} is $5$. $G=\Aut
(\Gamma)$ acts vertex- and edge-transitively on $\Gamma$.
For a vertex $x$ the action $G(x)^{\Gamma(x)}$ of the 
vertex stabilizer $G(x)$ on the set $\Gamma(x)$ of 
vertices adjacent to $x$ is doubly transitive
without regular normal subgroups. The kernel $G_1(x)$ of 
this action is
nontrivial.
\end{hypo}

In \cite{ivn2} and
\cite{ivn3} the classification  problem of graphs 
satisfying Hypothesis 6.1 was reduced 
to the classification of the
flag-transitive $P$-geometries. Now
we can formulate the final result (cf. \cite{ivn11}).
\begin{theorem}
Let $\Gamma$ be a graph satisfying Hypothesis \rom{6.1}. 
Then
$\Gamma$ is either the derived graph of one of the 
following 
$P$-geometries\,\rom{:} $\G(M_{22})$,
$\G(3 \cdot M_{22})$, $\G(Co_2)$, $\G(3^{23} \cdot Co_2)$, 
$\G(J_4)$,
$\G(BM)$, and $\G(3^{4371} \cdot BM)$, or it is a graph of 
valency $31$
which is related to the derived graph of $\G(J_4)$.
\end{theorem}

In the derived graph of $\G(M_{23})$ the condition $G_1(x) 
\ne 1$ fails.

The graphs in Theorem 6.2 are extremal by the order of the
vertex stabilizer in the general class of graphs with doubly
transitive action
$G(x)^{\Gamma(x)}$ \cite{tro}.
\begin{unique}
Construction of the geometries $\G(G)$ for
$G \cong J_4, BM$, and $M$ relies exclusively on the 
structure of the involution
centralizer $C=C_G(\tau)$ and on certain information on
 fusion in $G$ of involutions from
$O_2(C_G(\tau))$.
The characterization of the
geometries $\G(G)$ implies the following.
\end{unique}
\begin{theorem}
Let $G$ be a nonabelian simple group containing
an involution $\tau$ such that $C_G(\tau)$ is of the shape 
$2^{1+12}.3 \cdot
\Aut(M_{22})$, $2^{1+22}.Co_2$, or $2^{1+24}.Co_1$. 
Suppose that $C_G(O_2(C)) \le
O_2(C)$ and that $\tau^G \cap O_2(C) \ne \{\tau\}$. Then 
$G$ is uniquely
determined and is isomorphic to $J_4$, $BM$, or $M$, 
respectively.
\end{theorem}

For the original uniqueness proofs for the groups $J_4$, 
$BM$, and $M$ see
\cite{nor1,as1,ls,seg2,tho,nor2,gms}.

\begin{generators}
The classification of $T$- and $P$-geometries enables us
to obtain a characterization of certain sporadic simple 
groups stronger than
the characterization by the centralizer of an involution. 
The groups are
proved to coincide with the universal completions of 
certain of their subamalgams.
This provides us with presentations of the groups involved 
(the {\em geometric presentations}
\cite{ivn4}). In the case of $J_4$ the geometric 
presentation was proved in \cite{ivn7} to
be equivalent to a nice presentation conjectured by 
G.~Stroth and R.~Weiss in
\cite{sw}.
\end{generators}

In the case of $BM$ and $M$ the result establishes the 
correctness of the so-called
$Y$-presentations for these groups. The $Y$-presentations 
(cf. \cite{atlas,cns,nor3}) 
describe  groups as specific factor groups of Coxeter 
groups with diagrams 
having three
arms originating in a common node. The most famous is the 
$Y_{555}$-diagram with three 
arms of five edges each. The nodes on the arms are denoted 
by $a$, $b_i$, $c_i$, $d_i$, 
$e_i$, and $f_i$ for $i=$1, 2, and 3. After the 
announcement of
the geometric presentation of $M$ in \cite{ivn6}, 
S.~Norton \cite{nor4} proved its
equivalence to the corresponding $Y$-presentation. This 
resulted in the proof of
the following theorem conjectured by J.~Conway \cite{con2}.

\begin{theorem}
The Coxeter group corresponding to the $Y_{555}$-diagram 
subjected to a single additional relation 
$(ab_1c_1ab_2c_2ab_3c_3)^{10}=1$ is 
isomorphic to the wreath product $M \wreath 2$ of the 
Monster group and a group of order 
$2$ \rom{(}this wreath product is known as the 
Bimonster\rom{).}
\end{theorem}

The correctness of the $Y$-presentation for $BM$ is proved 
in \cite{ivn10}.

\medskip
\noindent
{\bf Representations and cohomologies}. Within the 
classification of
$T$- and $P$-geometries and their natural representations, 
considerable
information
on linear representations and nonsplit extensions of 
sporadic groups
was obtained. We formulate here only one result of this 
type which can be
deduced from \cite{is7}.
\begin{theorem}
Let $K$ be a field whose characteristic is not $2$.
Then $BM$ has a unique representation over $K$ of 
dimension $4371$. If 
$\CHAR(K)
\ne 3$, then the extension of $BM$ by the corresponding 
module always splits,
and for $K=GF(3)$ there is a unique nonsplit extension.
\end{theorem}

\noindent
{\bf Constructions}. Since the sporadic groups involved 
turned out to
coincide with the universal completions of certain of 
their subamalgams, there is
a possibility of producing an independent construction of 
these groups.
Such a construction is now in progress for $J_4$ where the 
corresponding
amalgam is realized in $GL_{1333}(C)$ \cite{im}. A 
possibility to use the
geometric characterization of the Monster group in order 
to simplify
its construction is discussed in \cite{ivn12}.
\section*{Acknowledgment}
We thank the referee for a number of helpful comments.


\begin{thebibliography}{[ABCDE]}

\bibitem[AS]{as1}
M.~Aschbacher and Y.~Segev, {\em The uniqueness of groups 
of type $J_4$\/}, Invent.
Math. {\bf 105} (1991), 589--607.

\bibitem[Bro]{bro}
K.S.~Brown, {\em Cohomology of groups\/}, Springer-Verlag, 
New York, 1982.

\bibitem[Bue1]{bue1}
F.~Buekenhout, Diagrams for geometries and groups, J. 
Combin. Theory Ser. 
A {\bf27} (1979), 121--151.

\bibitem[Bue2]{bue2}
F.~Buekenhout, ed., {\em Handbook on finite geometries \/} 
(to appear).

\bibitem[Con]{con2}
J.H.~Conway, $Y_{555}$ {\it and all that\/},
Groups,
Combinatorics and Geometry, Durham 1990 (M.~Liebeck and 
J.~Saxl, eds.),
London Math. Soc. Lecture Notes Ser., vol. {165},
 Cambridge
Univ. Press, 1992, pp. 22--23.

\bibitem[Atlas]{atlas}
J.H.~Conway, R.T.~Curtis, S.P.~Norton, R.A.~Parker and 
R.A.~Wilson,
{\em Atlas of finite groups\/}, Clarendon Press, Oxford, 
1985.

\bibitem[CNS]{cns}
J.H.~Conway, S.P.~Norton and L.H.~Soicher, {\it The 
bimonster, the group $Y_{555}$
and the projective plane of order\/} 3, Computers in 
Algebra (M. C.~Tangora,
ed.), Marcel Dekker, New York, 1988, pp. 27--50.

\bibitem[GMS]{gms}
R.L.~Griess, U.~Meierfrankenfeld, and Y.~Segev, {\em A 
uniqueness proof for the
Monster\/}, Ann. of Math. (2) {\bf 130} (1989), 567--602.

\bibitem[Hei]{hei}
S.~Heiss, {\em On a parabolic system of type\/} $M_{24}$, 
J. Algebra {\bf 142}
(1991), 188--200.

\bibitem[Ivn1]{ivn1}
A.A.~Ivanov, {\em On \rom{2}-transitive graphs of girth\/} 
5, European J. Combin. {\bf 8}
(1987), 393--420.

\bibitem[Ivn2]{ivn2}
\bysame, {\em Graphs of girth \rom{5} and diagram 
geometries related to the
Petersen graphs\/}, Soviet Math. Dokl. {\bf 36} (1988), 
83--87.

\bibitem[Ivn3]{ivn3}
\bysame, {\em The distance-transitive graphs admitting 
elations\/}, Math.
USSR-Izv. {\bf 35} (1990), 307--335.

\bibitem[Ivn4]{ivn4}
\bysame,  {\em Geometric presentation of groups with an 
application to the
Monster\/}, Proc. ICM-90, Kyoto, Japan, August 1990, 
Springer-Verlag,
New York,
1991, pp. 385--395.

\bibitem[Ivn5]{ivn5}
\bysame, {\em A geometric approach to the question of 
uniqueness for sporadic
simple groups\/}, Soviet Math. Dokl. {\bf 43} (1991), 
226--229.

\bibitem[Ivn6]{ivn6}
\bysame,  {\em A geometric characterization of the 
Monster\/},
Groups,
Combinatorics and Geometry, Durham 1990 (M.~Liebeck and 
J.~Saxl, eds.),
London Math. Soc. Lecture Notes Ser., vol. 165, Cambridge
Univ. Press., Cambridge, 1992, pp. 46--62.

\bibitem[Ivn7]{ivn7}
\bysame, {\em A presentation for $J_4$}, 
Proc. London Math. Soc. {\bf 64}
(1992), 369--396.

\bibitem[Ivn8]{ivn8}
\bysame, {\em A geometric characterization of 
Fischer\rom{'}s Baby Monster\/},
J. Algebraic Combin. {\bf 1} (1992), 43--65.

\bibitem[Ivn9]{ivn9}
\bysame,  {\em The minimal parabolic geometry of the 
Conway group $Co_1$
is simply connected\/}, Combinatorics 90: Recent Trends 
and Applications,
North-Holland, Amsterdam, 1992, pp. 259--273.

\bibitem[Ivn10]{ivn10}
\bysame, {\em Presenting the Baby Monster\/}, J. Algebra 
{\bf 163} (1994),
88--108.

\bibitem[Ivn11]{ivn11}
\bysame, {\em Graphs with projective subconstituents which 
contain
short cycles\/}, Surveys in Combinatorics (K.~Walker, 
ed.), London Math. Soc.
Lecture Notes Ser., vol. 187, Cambridge Univ. Press, 
Cambridge, 1993, 
pp. 173--190.

\bibitem[Ivn12]{ivn12}
\bysame, {\em Constructing the Monster via its 
$Y$-presentation\/},
Combinatorics, Paul Erd\"os is Eighty,  Bolyani Soc.
Math. Stud., vol. 1, Bolyani Math. Soc. Budapest, 1993, 
pp. 253--270.

\bibitem[IM]{im}
A.A.~Ivanov and U.~Meierfrankenfeld, 
{\em A construction and a uniqueness proof
for $J_4$}, preprint, 1994.

\bibitem[IPS]{ips}
A.A.~Ivanov, D.V.~Pasechnik, and S.V.~Shpectorov, {\em 
Non-abelian representations
of some sporadic geometries\/}, preprint, 1993.

\bibitem[IS1]{is1}
A.A.~Ivanov and S.V.~Shpectorov, {\em Geometries for 
sporadic groups related to
the Petersen graph.} I, Comm. Algebra {\bf 16} (1988), 
925--954.

\bibitem[IS2]{is2}
\bysame, {\em Geometries for sporadic groups related
to the Petersen graph.} II, European. J. Combin. {\bf 10} 
(1989), 347--362.

\bibitem[IS3]{is3}
\bysame, {\em The $P$-geometry for $M_{23}$ has no
nontrivial \rom{2}-coverings\/}, European J. Combin. {\bf 
11} (1990), 373--379.

\bibitem[IS4]{is4}
\bysame, {\em $P$-geometries of $J_4$-type have no
natural representations\/}, Bull. Math. Soc. Belg. S\'er. 
A {\bf 42} (1990),
547--560.

\bibitem[IS5]{is5}
\bysame, {\em Natural representations of the
$P$-geometries of $Co_2$-type\/}, J. Algebra 
{\bf164} (1994), 718--749.

\bibitem[IS6]{is6}
\bysame, {\em An infinite family of simply connected
flag-transitive tilde geometries\/}, Geom. Dedicata. {\bf 
45} (1993), 1--23.

\bibitem[IS7]{is7}
\bysame, {\em The last flag-transitive $P$-geometry\/},
Israel J. Math. {\bf 82} (1993), 341--362.

\bibitem[LS]{ls}
J.S.~Leon and C.C.~Sims, {\em The existence and uniqueness 
of a simple group
generated by \rom{3,4}-transpositions\/}, Bull. Amer. 
Math. Soc. {\bf 83} (1977),
1039--1040.

\bibitem[MS]{ms}
G.~Mason and S.D.~Smith, {\em Minimal \rom{2}-local 
geometries for the Held and
Rudvalis sporadic groups\/}, J. Algebra {\bf 79} (1982), 
286--306.

\bibitem[McK]{mck}
B.D.~McKay, Private communication, 1991.

\bibitem[Mei]{mei}
U.~Meierfrankenfeld, Private communication, 1991.

\bibitem[Nor1]{nor1}
S.P.~Norton, {\em The construction of\/} $J_4$,
The Santa Cruz Conference on Finite Groups
(B.~Cooperstein and G.~Mason, eds.),
Proc. Sympos. Pure Math., vol. 37,
Amer. Math. Soc., Providence, RI, 
 1980, pp. 271--278.

\bibitem[Nor2]{nor2}
\bysame, 
{\em The uniqueness of the Fischer-Griess Monster\/}, Finite
Groups---Coming of Age, Proc. 1982 Montreal Conf. 
(J.~McKay, ed.),
 Contemp. Math., vol.\ {45},
Amer. Math. Soc., Providence, RI, 1985, pp.\ 271--285.

\bibitem[Nor3]{nor3}
\bysame, {\em Presenting the Monster\,}?, Bull. Soc. Math. 
Belg.
S\'er. A {\bf 42} (1990), 595--605.

\bibitem[Nor4]{nor4}
\bysame, {\em Constructing the Monster\/}, Groups,
Combinatorics and Geometry, Durham 1990 (M.~Liebeck and 
J.~Saxl, eds.),
London Math. Soc. Lecture Note Ser., vol.\ {165}, 
Cambridge Univ. Press, Cambridge,
1992,  pp. 63--76.

\bibitem[Pas]{pas2}
A.~Pasini, {\em Geometries and diagrams\/}, Oxford Univ. 
Press,
London (to appear).

\bibitem[Ron1]{ron1}
M.A.~Ronan,  {\em Coverings of certain finite 
geometries\/}, Finite Geometries
and Designs, Cambridge Univ. Press, Cambridge, 1981,
pp. 316--331.

\bibitem[Ron2]{ron2}
\bysame, {\em Locally truncated buildings and $M_{24}$}, 
Math. Z. {\bf 180}
(1982), 489--501.

\bibitem[RSm1]{rsm1}
M.A.~Ronan and S.~Smith, 2-{\em local geometries for some 
sporadic groups\/}
(B.~Cooperstein and G.~Mason, eds.),
 Proc. Sympos. Pure Math., vol. 37,  Amer. Math. Soc., 
Providence,
RI, 1980, pp. 283--289.

\bibitem[RSm2]{rsm3}
\bysame, {\em Computation of \rom{2}-modular sheaves and
representations for $L_4(2)$, $A_7$, $3S_6$ and $M_{24}$}, 
Comm.
Algebra {\bf 17} (1989), 1199--1237.

\bibitem[RSt]{rst}
M.A.~Ronan and G.~Stroth, {\em Minimal parabolic 
geometries for the sporadic
groups\/}, European J. Combin. {\bf 5} (1984), 59--91.

\bibitem[Row1]{row1}
P.~Rowley, {\em On the minimal parabolic system related to 
$M_{24}$}, J. London
Math. Soc. {\bf 40} (1989), 40--56.

\bibitem[Row2]{row2}
\bysame, {\em Minimal parabolic systems with diagram
$\node\arc\node\arc\node\dstroke{\lower 
1.2ex\hbox{$\widetilde{}$}}\node$}, 
J. Algebra {\bf 141} (1991), 204--251.

\bibitem[Row3]{row3}
\bysame, {\em Pushing down minimal parabolic systems\/},
Groups,
Combinatorics and Geometry, Durham 1990 (M.~Liebeck and 
J.~Saxl, eds.),
London Math. Soc. Lecture Note Ser., vol. 165, Cambridge
Univ. Press, Cambridge,  1992, pp. 144--150.

\bibitem[Seg1]{seg1}
Y.~Segev, {\em On the uniqueness of the $Co_1$ 
\rom{2}-local geometry\/}, Geom.
Dedicata
{\bf 25}
(1988), 159--219.

\bibitem[Seg2]{seg2}
\bysame, {\em On the uniqueness of Fischer\rom{'}s Baby 
Monster\/}, Proc. London Math.
Soc. (3) {\bf 62} (1991), 509--536.

\bibitem[Shp1]{shp1}
S.V.~Shpectorov,  {\em A geometric characterization of the 
group $M_{22}$},
Investigations in Algebraic Theory of Combinatorial Objects,
Moscow, 1985, pp. 112--123; transl. by Kluwer Acad. Publ., 
1994.

\bibitem[Shp2]{shp2}
\bysame, {\em On geometries with diagram\/} $P^n$, 
preprint, 1988. 
(Russian)

\bibitem[Shp3]{shp3}
\bysame, {\em The universal \rom{2}-cover of the 
$P$-geometry ${\cal
G}(Co_2)$}, European J. Combin. {\bf 13} (1992), 291--312.

\bibitem[Shp4]{shp4}
\bysame, {\em Natural representations of some tilde and 
Petersen type
geometries\/}, Geom. Dedicata (to appear).

\bibitem[ShSt]{shst}
S.V.~Shpectorov and G.~Stroth, {\it Classification of 
certain types
of tilde geometries\/}, Geom. Dedicata {\bf49} (1994), 
155--172.

\bibitem[Smi]{smi}
S.D.~Smith, {\em Universality of the \rom{24}-dimensional 
embedding of the
\rom{.1 2}-local
geometry\/}, Comm.  Algebra (to appear).

\bibitem[Str]{str}
G.~Stroth, {\em Parabolics in finite groups\/}, Proc. 
Rutgers Group Theory
Year, 1983--1984, Cambridge Univ. Press, Cambridge, 1984, 
pp. 211--224.

\bibitem[SW]{sw}
G.~Stroth and R.~Weiss, {\em Modified Steinberg relations 
for the group\/} $J_4$,
Geom. Dedicata {\bf 25} (1988), 513--525.

\bibitem[Tho]{tho}
J.G.~Thompson, {\em Uniqueness of the Fischer-Griess 
Monster\/}, Bull. London Math.
Soc. {\bf 11} (1979), 340--346.

\bibitem[Tim]{tim}
F.~Timmesfeld, {\em Classical locally finite Tits chamber 
systems of rank\/} 3, J.
Algebra {\bf 124} (1989), 9--59.

\bibitem[Tit1]{tit1}
J.~Tits, {\em Buildings of spherical type and finite 
BN-pairs\/}, Lecture Notes
in Math., vol.\ {386}, Springer-Verlag, New York, 1974.

\bibitem[Tit2]{tit2}
\bysame, {\it A local approach to buildings\/}, 
The Geometric Vein (Coxeter-%
Festschrift), Springer-Verlag, Berlin, 1982, pp. 519--547.

\bibitem[Tro]{tro}
V.I.~Trofimov, {\em Stabilizers of the vertices of graphs 
with projective
suborbits\/}, Soviet Math. Dokl. {\bf 42} (1991), 825--827.


\end{thebibliography}
\end{document}